# DERANGEMENTS IN SYMMETRIC COST MATRICES

## Howard Kleiman

**I. Introduction.** Let M be an $n \times n$ symmetric cost matrix. Assume that $D$ is a derangement of edges in M, i.e., a set of disjoint cycles consisting of edges that contain all of the $n$ points of M. The modified Floyd-Warshall algorithm applied to $(D')^{-1} A^-$ (where $A$ is an asymmetric cost matrix containing $D'$, a derangement) yielded a solution to the Assignment Problem in $O(n^2 \log n)$ running time. Here, applying a variation of the modified F-W algorithm to $D^{-1} M^-$, we can obtain $D = D_{FWABS}$, the smallest derangement obtainable using the modified F-W algorithm. Let $T_{TSPOPT}$ be an optimal tour in M. We also give conditions for obtaining $D_{ABSOLUTE}$, the smallest derangement obtainable in M, where $|D_{ABSOLUTE}| \leq |T_{TSPOPT}|$.

**II. Theorems.**

**Theorem 1.** *Let M be a value matrix containing both positive and negative values. Suppose that M contains one or more negative cycles. Then if a negative path P becomes a non-simple path containing a negative cycle, $C$, as a subpath, C is obtainable as an independent cycle in the modified F-W algorithm using fewer columns than the number used by P to construct $C$.*

**Proof.** Using the modified F-W algorithm, no matter which vertex, $a_i$, of $C$ is first obtained in $P$, it requires the same number of columns to return to $a_i$ and make $P$ a non-simple path. Since $b_1 \neq a_1$, it requires at least one column to go from $b_i$ to $a_1$, concluding the proof.

In the following example, we construct $P$ and $C$.

**Example 1**

$P = [1^{-20} 3^5 7^{-5} 13^{12} 15^1 19^3 20^{-18} 18^1 14^1 6^3 7^{-5}]$,

$C = (20\ 18\ 14\ 6\ 7\ 13\ 15\ 19)$.

In what follows, Roman numerals represent numbers of iterations.

| P | C |
|---|---|
| **I.**   [1 3] + 3 | **I.**   [20 18] + 18 |
|      [1 7] + 4 | ------------------------ |
|      [1 13] + 6 | **II.**  [20 14] + 16 |
|      [1 15] + 2 | ------------------------ |
|      [1 19] + 4 | **III.** [20 6] + 12 |
|      [1 20] + 1 |      [20 7] + 1 |
| ------------------------ |      [20 13] + 6 |
| **II.**  [1 18] + 18 |      [20 15] + 2 |
| ------------------------ |      [20 19] + 4 |
| **III.** [1 14] + 16 |      [20 20] + 1 |
| ------------------------ | ------------------------ |
| **IV.** [1 6] + 12 | **+ 60** |
|      [1 7] + 1 | |
| ------------------------ | |
| **+ 67** | |

In the next theorem, $P_{in}$ is an $n \times n$ matrix, that gives all paths obtained while we are in the i-th iteration of the modified F-W algorithm. *NEGVALUES* is an $n \times 2$ matrix satisfying the following conditions:

(a) Its first column consists of the numbers from 1 through n.

(b) If $C$ is a negatively-valued cycle, then its value is placed in the second column of each row number that is one of its points.

As each point is added to a row, the values in the second column are sorted in increasing order of magnitude.

$|P|$ represents the sum of the values in $\sigma^{-1}M^{-}(in)$ of the arcs of path $P$.

**Theorem 2.** Let $D$ be a derangement of edges in M, while $Q$ is a simple path of permutation arcs in $P_{in}$ obtained by applying the modified F-W to $D^{-1}M^{-}$. If $Q = (d \ ... \ c')$ and we wish to consider adjoining the arc $(c' \ a)$ to $Q$, suppose that

    (1)

(i) $D(a)$ is not the $D^{-1}(c')-th$ entry in row $d$, or

(ii) if it is, then $(D(a)\ D^{-1}(c'))$ is not an arc of $Q$.

(2) any of the following holds: (i) the entry $(d\ a)$ of $P_{in}$ is blank; (ii) if the entry $(d\ a)$ exists in $P_{in}$, and $|(d\ a)|$ is the value of $(d\ a)$ in $\sigma^{-1}M^{-}(in)$, then $|(d\ ...\ c'\ a)|-|(d\ a)|$ is not a number in the second column of the $a-th$ row of $NEGVALUES$; (iii) $a$ isn't a point of $Q$.

Then $Q' = (d\ ...\ c'\ a)$ is a simple path of permutation arcs.

**Proof.** Assuming our conclusion is true for $Q$, we must prove that $(c'\ a)$ yields a directed edge not symmetric to a directed edge obtained by an earlier arc of $Q$ and that $Q'$ doesn't contain a cycle.

(1) (i) and (ii) guarantee that even if an edge symmetric to $(c'\ D(a))$ can be obtained from an arc in row $d$, it doesn't belong to $Q$ and, therefore, to $Q'$.

(2) Using modified F-W, we can't obtain a non-negative path. Assuming that $Q$ doesn't contain a negative cycle, it can only occur in $Q'$ if $a$ is a point of $Q$. If that is the case and the entry $(d\ a)$ of $P_{in}$ isn't blank, then a negative cycle, $C$, exists starting and ending at $(d\ a)$. $|C|=|(d\ ...\ c'\ a)|-|(d\ a)|$. Thus, $|C|$ would lie in the second column of the $a-th$ row of $NEGVALUES$. Given the validity of condition 2(ii), it is possible that 2(iii) is the case, i.e., $Q$ doesn't contain $a$. We ascertain this by back-tracking in $P_{in}$ from $c'$ to $d$.

**Theorem 3.** Let $S$ be a set of pair-wise point-disjoint cycles constructed from paths satisfying theorem 2. Let $SYM$ be the subset of arcs of form $(a_i\ b_i)$ in the cycles of $S$ each of which is symmetric to an arc of $D$. Assume that the following are valid:

(1) for each subscript $i$, $D(b_i)$ is a point of a cycle in $S$;

(2) for each of the arcs $(d_j\ e_j)$ of cycle $C_j$ in $S - SYM$, conditions (i) and (ii) of theorem 2 hold with respect to all of the arcs in $S - SYM - C_j$.

Then if $s$ is the permutation that is the product of the cycles in $S$, $Ds$ is a derangement of edges.

**Proof.** From (1), we know that given an arc $(a_i\ b_i)$ in $SYM$, there must exist two edges - one obtainable from an arc in $SYM$, the other from an arc in $S$ - that are of the form $[a_i\ D(b_i)][D(b_i)\ c_i]$, $c_i \neq a_i$. Thus, we get two distinct edges when $D$ is applied to $(a_i\ b_i)$ and $(D(b_i)\ D^{-1}(c_i))$, respectively. If $(d_j\ e_j)$ is an arc of $C_j$, we know that all of the remaining arcs of $C_j$ yield different edges than that obtained from $(d_j\ e_j)$. We thus need only be concerned with arcs lying on other cycles that don't lie in $SYM$. Thus, by applying (i) and (ii) to all of the arcs in $S - SYM - C_j$, we guarantee that no two arcs yield the same edge. Since $s$ is a permutation, $Ds$ is also a permutation. By construction, applying $D$ to any product of disjoint cycles in $D^{-1}M^-$ always yields a derangement if we include the possibility of 2-cycles. Theorem 2 and (1) and (2) of theorem 3 eliminate the occurrence of 2-cycles (edges). It follows that $Ds$ is a derangement of edges.

**Corollary 1.** Let $T$ be a tour that is an upper bound for $T_{FWTSPOPT}$. Then the following holds:

$$|D_{FWABS}| \leq |T_{FWTSPOPT}| \leq |T_{UPPERBOUND}|$$

**Proof.** $T_{FWTSPOPT}$ is a derangement of edges obtained using the modified F-W algorithm. Therefore, it cannot have a smaller value than $D_{FWABS}$.

*Note 1.* Before beginning the construction of paths, we write both $D$ and $D^{-1}$ in row form, i.e., the top row is 1, 2, ... , n; the bottom row is either the corresponding value of $D$ or $D^{-1}$. Using them, we can quickly obtain the values $D(a)$ and $D^{-1}(b)$.

*Note 2..* Suppose we obtain the negatively-valued permutation $s_1$ in $D^{-1}M^-$ where $s_1$ satisfies the conditions of theorem 3. Then $Ds_1 = D_1$ has the property that $|D_1| < |D|$. We put those edges containing the initial points of $s_1$ in row form. We then construct $D_1^{-1}M^-$. If we can obtain a negatively-valued permutation, $s_2$, in $D_1^{-1}M^-$, $D_1s_2 = D_2$. We then construct $D_2^{-1}M^-$. We continue until we reach $D_m$ where $D_m^{-1}M$ contains no negatively-valued permutations satisfying theorem 3. In this case, $D_m = D_{FWABS}$.

**Theorem 4.** Let $S_i$, $i = 1,2, ... ,r$ be all sets of pair-wise point-disjoint negatively-valued cycles obtainable in $D_{FWABS}^{-1}M^-$. Let $V = min\ \{|S_i|\ |\ i = 1,2,... ,r\}$. Define $POS$ to be the set of all sets $P_j$, $j = 1, 2, \ldots, s$ of positive cycles in $D_{FWABS}^{-1} M^-$ of value less than $-V$ that are pair-wise point-disjoint to each other as well as to the cycles of at least one $S_i$. Assume each such pair satisfies the conditions of theorem 3. Let $S_i^* \bigcup P_j^*$ be the set of permutations of smallest value among all such pairs of sets, while $sp$ is the product of its cycles. Then $D_{FWABS} sp = D_{ABSOLUTE}$.

**Proof.** Using F-W, we can obtain all negatively-valued cycles in $D_{FWABS}^{-1} M^-$. To obtain the positively-valued cycles, we could keep track of all paths that we obtain while searching for cycles of value less than $-V$. In particular, before we delete an entry in $P_{in}$, we store the path that is being replaced. When we have obtained all positively-valued cycles of value less than $-V$ using F-W, we first delete all paths whose initial point is not a determining vertex of a cycle. We then try to extend the remaining paths (using each determining vertex initial point as a root of a tree) to see if we can obtain further cycles of value less than $-V$. Once we have obtained all possible positively-valued cycles, we see if we can obtain a negatively-valued permutation satisfying the conditions of theorem 3. If more than one exists, we choose that which has the smallest negative value. Calling this permutation $sp$, $D_{FWABS} sp = D_{ABSOLUTE}$.